\documentclass{article}
\usepackage[arrow, matrix, curve]{xy}
\usepackage{amsthm}
\usepackage{hyperref}
\input{amssym.tex}
\newtheorem{thm}{Theorem}
\newtheorem{lem}{Lemma}
\newtheorem{cor} {Corollary}
\newtheorem{exm} {Example}

\newtheorem{pro} {Proposition}
\newtheorem{df} {Definition}

	\newenvironment{pf}{{\it Proof:}\quad}{\hfill$QED$}

\begin{document}

	\title{Invariants preserved by mutation} 
	\author{Thilo Kuessner}
	\date{}
	\maketitle
	\begin{abstract} \noindent
 We prove that generalized mutation preserves several geometric invariants such as the volume 
and Goncharov invariant of closed or Q-rank 1 locally symmetric spaces.
\noindent
		\end{abstract}

Rigidity results in the theory of locally symmetric spaces imply that geometrically defined invariants as the volume are topological invariants. Yet it remains largely mysterious how these invariants are determined by the topology and how they behave with respect to topological operations such as cut and paste.  Hereby "cut and paste" refers to the following operation: we have a  properly embedded, 2-sided, codimension 1
submanifold $\Sigma$ in a compact manifold $M$ and we denote by $M^\tau$ the manifold which is
the result of cutting $M$ along $\Sigma$ and regluing via a diffeomeorphism $\tau:\Sigma\rightarrow\Sigma$.

In \cite{rub} Ruberman considered the case of hyperbolic 3-manifolds and showed that mutation of a hyperbolic link yields a hyperbolic link of the same volume. More generally he proved that for hyperbolic $3$-manifolds $M$ and certain pairs $\left(\Sigma,\tau\right)$, especially 
for the hyperelliptic involution of the genus 2 surface, $M^\tau$ is always hyperbolic with $vol\left(M^\tau\right)=vol\left(M\right)$ if $\Sigma\subset M$ is incompressible and boundary-incompressible. (The latter conditions are needed only to guarantee hyperbolicity of $M^\tau$.) Neumann indicated in \cite{neu} that also the $PSL(2,{\bf C})$-fundamental class and up to torsion the Bloch invariant of hyperbolic $3$-manifolds is preserved under mutation. In \cite{ku2} we gave a topological proof of Ruberman's theorem using the fundamental class construction. In \cite{kkk} we used an analogous argument to prove that also the volume of flag structures is preserved under mutation. The aim of this paper is to prove in a general setting that $G$-fundamental classes and hence various geometric invariants are preserved under generalized mutation.

For a closed, orientable $d$-manifold $M$ and a representation $\rho:\pi_1M\rightarrow G$ one has the naturally associated $G$-fundamental class $(B\rho)_*\left[M\right]\in H_d(BG)$. If $M$ has $\pi_1$-injective boundary, $G$ is a semisimple Lie group without compact factor, and $\rho$ sends $\pi_1\partial_iM$ to a parabolic subgroup $P_i\subset G$ for each component $\partial_iM$ of $\partial M$, then we can still associate a fundamental class $(B\rho)_*\left[M,\partial M\right]\in H_d(BG^{comp})$ for a certain completion $BG^{comp}$, see Section 3.1. We discuss in Section 2 that this fundamental class determines several geometric invariants. The generalized mutations to be considered will be defined in \hyperref[genmut]{Definition \ref*{genmut}}. The general result proved in this paper is the following theorem.

\begin{thm}\label{Thm1} Assume $M$  is a compact, oriented $d$-manifold such that $int(M)$ is a ${\bf Q}$-rank $1$ locally symmetric space of noncompact type.
Let $\rho:\pi_1M\rightarrow G$ be a representation to a semisimple Lie group without compact factor, sending the fundamental group of each boundary component $\pi_1\partial_iM$ to some parabolic subgroup $P_i\subset G$.

If $\rho^\tau:\pi_1M^\tau\rightarrow G$ is a generalized mutation of $\rho$, then 
$$(B\rho)_*\left[M,\partial M\right]_{\bf Q}=(B\rho^\tau)_*\left[M^\tau,\partial M^\tau\right]_{\bf Q}\in H_d(BG^{comp};{\bf Q}).$$\end{thm}
\noindent
The corresponding result for ${\bf Z}$-coefficients is not true, examples are given in \cite{neu}.
The following invariants are determined by the rational $G$-fundamental class, see Section 2.
\begin{cor}\label{cor}The volume, the (generalized) rational Bloch invariants and Goncharov invariants of ${\bf Q}$-rank 1 locally symmetric spaces, as well as the rational Bloch invariant and volume of CR structures and flag structures are preserved under generalized mutations.\end{cor}

\section{Generalized mutation}

\subsection{Definition}
In this paper we will consider the following situation. We have a compact manifold $M$ (possibly with boundary) and 
a properly embedded, 2-sided, codimension $1$
submanifold $\Sigma\subset M$. 
We consider a diffeomorphism $\tau:\left(\Sigma,
\partial\Sigma\right)\rightarrow \left(\Sigma,\partial\Sigma\right)$ and we denote $M^\tau$ the manifold which is
the result of cutting $M$ along $\Sigma$ and regluing via $\tau$.

We will also assume that a representation $\rho:\pi_1M\rightarrow G$ is given. This representation may arise from the identification of $\pi_1M$ with a discrete subgroup $\Gamma\subset G$ in case that $M=\Gamma\backslash G/K$ is a locally symmetric space of noncompact type, but we will also be interested in other representations, e.g. arising from CR or flag structures.

\begin{df}\label{genmut}({\bf Generalized mutation}) For a fixed representation $\rho:\pi_1M\rightarrow G$ we say that a finite order homeomorphism $\tau:\Sigma\rightarrow\Sigma$  is a generalized mutation if there exists some finite order $A\in G$ with $$\rho\left(\tau_*h\right)=A\rho\left(h\right)A^{-1}$$
for all $h\in\pi_1\Sigma$.\end{df}

We remark that the condition of $A$ having finite order is implied 
if $\rho$ is the holonomy of a hyperbolic $3$-manifold (\cite[Observation 1.2]{ku2}) or if $\rho(\pi_1\Sigma)$ is Zariski-dense.

\begin{lem}\label{X} If $\tau:\Sigma\rightarrow\Sigma$  is a generalized mutation for a representation $\rho:\pi_1M\rightarrow G$, then there is a representation
$$\rho^\tau:\pi_1M^\tau\rightarrow G$$
such that the retrictions of $\rho$ and $\rho^\tau$ to $\pi_1(M - \Sigma)$ agree.\end{lem}
\begin{pf} We can choose regular neighborhoods of $\Sigma$ in $M$ and $M^\tau$ and an identification of their complements. We denote $X$ the union of $M$ and $M^\tau$ along this identification. 
Let $\pi_1M=<S\mid R>$ be a presentation of $\pi_1M$. The Seifert-van Kampen Theorem implies
$$\pi_1X=<S,t\mid R, tht^{-1}=\tau_*\left(h\right)\ \forall\ h\in\pi_1\Sigma>.$$
\hyperref[genmut]{Definition \ref*{genmut}} implies that we can extend $\rho$ to a representation $\rho^X:\pi_1X\rightarrow G$ by defining $\rho(t)=A$. Composition with $\pi_1M^\tau\rightarrow\pi_1X$ yields $\rho^\tau:\pi_1M^\tau\rightarrow G$.\end{pf}

\subsection{Examples}
{\bf Symmetries of Representation Variety.} The mapping class of $\tau:\Sigma\rightarrow\Sigma$ acts on $R(\pi_1\Sigma,G)=Hom\left(\pi_1\Sigma,G\right)/\sim$, the space of representations 
up to conjugation. The condition from \hyperref[genmut]{Definition \ref*{genmut}} is obviously equivalent to the conjugacy class of $\rho\mid_{\pi_1\Sigma}$ being a fixed point for $\tau_*$ on $R(\pi_1\Sigma,G)$. We do not know a general approach for finding such fixed points but a remarkably general case is provided by the following example due to Ruberman. (Here $R_{d,f}(\pi_1\Sigma,G)$ is the subset of discrete, faithful, parabolics-preserving representations upon conjugation, an analogue of Teichm\"uller space.)

\begin{exm}\label{ruberman}{\bf Surfaces in hyperbolic 3-manifolds.} 
If $\Sigma\subset M$ is the genus $2$ surface and $\tau:\Sigma\rightarrow \Sigma$ the hyperelliptic involution, then by \cite[Theorem 2.2]{rub} $\tau_*$ acts trivially on $R_{d,f}(\pi_1\Sigma,SL(2,{\bf C}))$, every $\rho$ is a fixed point. The same is true for the $\tau$-invariant subsurfaces of $\Sigma$: the 1- and 2-punctured torus and the 3- and 4-punctured sphere.
\end{exm}

\begin{exm}\label{total}{\bf Totally geodesic submanifolds.}
If $M=\Gamma\backslash G/K$ is a closed locally symmetric space of noncompact type of dimension $d\ge 4$, $\Sigma\subset M$ a totally geodesic hypersurface and $\tau:\Sigma\rightarrow \Sigma$ a diffeomorphism, then 
by \cite[Theorem IV.7.1]{hel} $\Sigma$ is a locally symmetric space, upon conjugation we can assume $\Sigma=H/K$ for some subgroup $H\subset G$, hence Mostow rigidity means that $\tau:\Sigma\rightarrow  \Sigma$ is homotopic to an isometry (of finite order) and there is some $A\in H\subset G$ with $\tau_*\left(h\right)=AhA^{-1}$
for all $h\in \pi_1\Sigma$.\end{exm}

\subsection{Discreteness}

Some invariants considered in this paper are defined only for discrete embeddings into a Lie group (this is the case for the volume or the Bloch invariant), some are defined for other representations (this is the case for the volume of flag structures or the Goncharov invariant). For the first case it is natural to ask whether the mutation $\rho^\tau$ of a discrete embedding $\rho:\pi_1M\rightarrow G$ is again discrete. In \cite{rub} Ruberman answered this positively for lattices in $SL(2,{\bf C})$. When $\rho(\pi_1\Sigma)$ is geometrically finite (hence $1$-quasifuchsian, \cite[Definition 9.2]{ka}), then in \cite[Proposition 3.1]{ku2} we gave another proof by using the Maskit combination theorem from \cite[Chapter VII]{mas}. The Maskit combination theorem in the formulation of \cite{mas} has an exact generalisation to higher-dimensional hyperbolic manifolds by the recent work of Li-Ohshika-Wang (\cite[Theorem 4.2]{low}). Thus one can literally adapt the proof of \cite[Proposition 3.1]{ku2} to obtain the following result:
\begin{pro} If $M$ is a compact, oriented $n$-manifold such that $int(M)$ is hyperbolic with holonomy $\rho:\pi_1M\rightarrow Isom^+(H^n)$, and if $\Sigma$ is a properly embedded, 2-sided, codimension 1 submanifold such that $\rho(\pi_1\Sigma)$ is geometrically finite, $(n-2)$-quasifuchsian, then $\rho^\tau(\pi_1M^\tau)$ has discrete image for any generalized mutation $\tau:\Sigma\rightarrow\Sigma$.
\end{pro}

\subsection{The closed case}
Though this is superfluous from a logical point of view - as we handle the more general ${\bf Q}$-rank 1 case 
in Section 3 - we first discuss the, technically simpler, proof of \hyperref[Thm1]{Theorem \ref*{Thm1}} for the closed case. The idea is essentially due to Neumann (\cite{neu}, in the context of hyperbolic $3$-manifolds) and in this case the simplicity and beauty of the argument can hopefully be better appreciated.

Let $X$ be constructed as in the proof of \hyperref[X]{Lemma \ref*{X}}. The images of the fundamental classes of $M,M^\tau$ and the mapping torus $T^\tau$ satisfy the relation

$$\left[M\right]-\left[M^\tau\right]= \left[T^\tau\right]\in H_d\left(X\right).$$
From the proof of \hyperref[X]{Lemma \ref*{X}} we have a representation $\rho^X:\pi_1X\rightarrow G$ with $\rho^X\mid_{\pi_1M}=\rho$ and $\rho^X\mid_{\pi_1M^\tau}=\rho^\tau$. Then
%
%
$$(B\rho)_*\left[M\right]-(B\rho^\tau)_*\left[M^\tau\right]=(B\rho^X)_*\left[T^\tau\right]$$
and thus Theorem 1 will follow once we have proved $(B\rho^X)_*\left[T^\tau\right]_{\bf Q}=0$.
%
%

$\tau$ and $A$ have finite order, say $\tau^n=id$ and $A^n={\bf 1}$ for an $n\in{\bf N}$. Hence we have an $n$-fold covering $p_n:\Sigma\times S^1\rightarrow T^\tau$.
In the following diagram we denote by $i_\Sigma$ and $i_{T^\tau}$ the inclusions and by $P_1$ the projection to the 1st factor.
$$\begin{xy}\xymatrix{ \Sigma\times S^1\ar[d]^{p_n}\ar[r]^{P_1}&\Sigma\ar[d]^{i_\Sigma}\\
T^\tau\ar[r]^{i_{T^\tau}}&
X\ar[r]^{B\rho^X}&
\mid BG\mid}
\end{xy}$$

The left-hand square does not commute, even homologically, but we claim that the compositions with $B\rho^X$ commute homologically. 

Since $\mid BG\mid$ is aspherical it suffices to look at the fundamental group.  
We have $\pi_1\left(\Sigma\times S^1\right)=\pi_1\Sigma\oplus{\bf Z}$ and for the $\pi_1\Sigma$-summand of course already the left-hand square is commuting.
Moreover, for the generator $s\in \pi_1S^1={\bf Z}$ we have $P_{1*}\left(s\right)=0$, on the other hand $(i_{T^\tau}p_{n})_*\left(s\right)=t^n$ and thus
$\left(B\rho_X i_{T^\tau}p_n\right)_*  \left(s\right)=\rho_X\left(t^n\right)=A^n={\bf 1}$.
Therefore the induced homomorphisms of fundamental groups commute, the same is trivially true for higher homotopy groups, and Hurewicz' Theorem implies commutativity in homology.

With $d=dim(M)$ we have $H_d(\Sigma)=0$. The homomorphism $\left( B\rho^Xi_\Sigma P_1\right)_d$ factors over
$H_d\left(\Sigma\right)$ and is therefore trivial. By the discussion before this implies
$$\left( B\rho^X
i_{T^\tau}p_n\right)_*\left[\Sigma\times {\bf S}^1\right]=0.$$
We have $p_{n*}\left[\Sigma\times S^1\right]=n\left[T^\tau\right]$, thus {\em rationally} 
$\left( B\rho^X \right)_*\left[T^\tau\right]_{\bf Q}
=0,$ hence the claim.

\section{Invariants obtained from the $G$-fundamental class}

For a manifold $M$ and a group $G$, 
a representation $\rho:\pi_1M\rightarrow G$ induces a continuous map
$B\rho:M\rightarrow \mid BG\mid$. If $M$ is a closed, oriented $d$-manifold with fundamental class $\left[M\right]\in H_d(M)$, then the $G$-fundamental class of $(M,\rho)$ is defined as $(B\rho)_*\left[M\right]\in H_d(BG)$.

We will be more generally interested in the situation that $M$ has ($\pi_1$-injective, possibly disconnected) boundary, $G$ is a semisimple Lie group without compact factors and that the representation $\rho:\pi_1M\rightarrow G$ maps the fundamental group of each component $\partial_iM$ to some parabolic subgroup $P_i\subset G$. 
Then we have an element 
$$(B\rho)_*\left[M,\partial M\right]\in H_d(BG^{comp}),$$
see Section 3.1 for its construction.
Several invariants can be derived from this element.

{\bf Volume of locally symmetric spaces.} If $\rho:\Gamma\rightarrow G$ is the inclusion of a ${\bf Q}$-rank $1$ lattice and $int(M)=\Gamma\backslash G/K$ the locally symmetric space, then by \cite[Section 4.2.3]{ku} there is the extended volume cocycle $\overline{c\nu}_d$ defined by $\overline{c\nu}_d(g_1,\ldots,g_d)=\int_{str(\tilde{x},g_1\tilde{x},\ldots,g_1\ldots g_d\tilde{x})}dvol$ for $(g_1,\ldots,g_d)\in BG$ and by $\overline{c\nu}_d(p_1,\ldots,p_{d-1},c_i)=\int_{str(\tilde{x},p_1\tilde{x},\ldots,p_1\ldots p_{n-1}\tilde{x},c_i)}dvol$ for $(p_1,\ldots,p_{n-1},c_i)\in Cone(BG)$. Its cohomology class does not depend on $\tilde{x}\in G/K$ and by \cite[Lemma 6.3]{kkk} we have
$$<\left[\overline{c\nu}_d\right],(B\rho)_*\left[M,\partial M\right]>=vol(M).$$
In particular the volume is determined by $(B\rho)_*\left[M,\partial M\right]_{\bf Q}$,

{\bf Goncharov invariant.} If again $int(M)=\Gamma\backslash G/K$ is a ${\bf Q}$-rank one locally symmetric space of finite volume, then by Weil rigidity we can assume $\Gamma\subset G(\overline{\bf Q})$. If $\rho:\Gamma\rightarrow SL(N,\overline{\bf Q})$ is the restriction of some representation $G(\overline{\bf Q})\rightarrow SL(N,\overline{\bf Q})$, then by \cite[Proposition 7.1]{kkk} the rational fundamental class $(B\rho)_*\left[M,\partial M\right]_{\bf Q}\in H_d(BSL(N,\overline{\bf Q})^{comp},{\bf Q})$ has a (unique) preimage $\overline{\gamma}(M)\in  H_d(BSL(N,\overline{\bf Q}),{\bf Q})$ and also an associated element $\gamma(M)\in K_d(\overline{\bf Q})\otimes{\bf Q}$, which are called the homological and K-theoretic Goncharov invariant, respectively.
Similarly, if ${\bf F}\subset{\bf C}$ is a subring with 1 and $\rho(\Gamma)\subset SL\left(N,{\bf F}\right)$, then one obtains $\gamma(M)\in K_d({\bf F})\otimes{\bf Q}$. By definition, $\gamma(M)$ is determined by $(B\rho)_*\left[M,\partial M\right]_{\bf Q}$.

{\bf Bloch invariant.} It is proved in the appendix of \cite{ds} that the Bloch-Wigner morphism sends the $PSL(2,{\bf C})$-fundamental class of a closed hyperbolic $3$-manifold to its Bloch invariant $\beta(M)$. The corresponding result for cusped hyperbolic $3$-manifolds as well as for the generalized Bloch invariant of (closed or ${\bf R}$-rank $1$) locally symmetric spaces has been proved in \cite[Theorem 1]{ku3}, see also \cite[Definition 8.1]{kkk} for ${\bf Q}$-rank $1$ spaces. Consequently $(B\rho)_*\left[M,\partial M\right]_{\bf Q}$ is mapped to the {\em rational Bloch invariant} $\beta(M)\otimes 1$. 

{\bf Bloch invariant of CR structures (\cite{FW}).} If $M$ is a finite-volume hyperbolic manifold and $\rho:\pi_1M\rightarrow SU(2,1)$ a reductive representation, then \cite[Lemma 10.3]{kkk} implies that the Bloch invariant $\beta_{FW}(M)$ is determined by $(B\rho)_*\left[M,\partial M\right]$.

{\bf Bloch invariant and volume of flag structures (\cite{bfg}).} If $M$ is a finite-volume hyperbolic $3$-manifold, and $h:{\bf C}P^1\rightarrow {\mathcal{F}}l({\bf C}^3)$ is equivariant with respect to some representation $\rho:\pi_1M\rightarrow SL(3,{\bf C})$, then Bergeron-Falbel-Guilloux define a Bloch invariant $\beta_h(M)$ which generalizes the Bloch invariants of hyperbolic and CR structures. It follows from \cite[Lemma 10.7]{kkk} that $\beta_h(M)$ is determined by $(B\rho)_*\left[M,\partial M\right]$. In particular
the volume of the flag structure is determined by $(B\rho)_*\left[M,\partial M\right]_{\bf Q}$.

\section{Proof}
\subsection{Recollections from \cite[Section 4]{ku}: Construction of $(B\rho)_*\left[M,\partial M\right]$}
\begin{df}\label{dc}
For a manifold $M$ and $\partial_1M,\ldots,\partial_sM$ the components of $\partial M$,
we let
$$DCone\left(\cup_{i=1}^s\partial_iM\rightarrow M\right)$$ 
be the union along $\partial M$ of $M$ with the (disjoint) cones over $\partial_1 M,\ldots,\partial_s M$.

For a group $\Gamma$ and a collection of subgroups $\Gamma_1,\ldots,\Gamma_s$, we denote by
$$B\Gamma^{comp}:=DCone\left(\cup_{i=1}^sB\Gamma_i\rightarrow B\Gamma\right)$$ 
the union along $\cup_{i=1}^sB\Gamma_i$
of the simplicial set $B\Gamma$ with the cones over $B\Gamma_i$.

If $G$ is a semisimple Lie group without compact factor and $\partial_\infty G/K$ the Hadamard boundary of the associated symmetric space, then we denote $$BG^{comp}=
DCone\left(\dot{\cup}_{c\in\partial_\infty G/K}BG\rightarrow BG\right).$$ \end{df}

Let $M$ be connected and the components
$\partial_1M,\ldots,\partial_sM$ be $\pi_1$-injective. We fix $x\in M,
x_i\in\partial_iM$ and pathes $l_i:\left[0,1\right]\rightarrow M$ witn $l_i\left(0\right)=x, l_i\left(1\right)=x_i$ for $1\le i\le s$. Conjugation by $l_i$ identifies $\pi_1\left(\partial_iM,x_i\right)$ 
to a subgroup $\Gamma_i\subset \Gamma:= \pi_1\left(M,x\right)$. 

Let $\mid B\Gamma
\mid$ denote the geometric realization of the simplicial set $B\Gamma$. The classifying map $\Psi_M:M\rightarrow \mid B\Gamma\mid$ extends to $\Psi_M:DCone\left(\cup_{i=1}^s\partial_iM\rightarrow M\right)\rightarrow B\Gamma^{comp}$. 

If $int\left(M\right)\cong\Gamma\backslash G/K$ is a ${\bf Q}$-rank $1$ locally symmetric space, then by \cite[Lemma 6.2]{kkk} composition of $\left(\Psi_M\right)_*$ with the isomorphism $H_*\left( \mid
B\Gamma^{comp} \mid\right)\cong H_*^{simp}\left(
B\Gamma^{comp}\right)$ yields the inverse of the Eilenberg-MacLane-isomorphism:
$$EM^{-1}:H_*\left(DCone\left(\cup_{i=1}^s\partial_iM\rightarrow M\right)\right)\rightarrow H_*^{simp}
\left(B\Gamma^{comp}\right).$$

For $*\ge 2$ we have
an isomorphism $H_*\left(M,\partial M\right)\cong H_*\left(DCone\left(\cup_{i=1}^s\partial_iM\rightarrow M\right)
\right)$, see \cite[Lemma 5.1]{kkk}. Thus, if $dim\left(M\right)\ge 2$, then
$EM^{-1}\left[M,\partial M\right]\in H_*^{simp}\left(B\Gamma^{comp}\right)$ 
is defined. 

We assume that $\rho:\Gamma\rightarrow G$ is a representation to a simple Lie group without compact factor, sending each $\Gamma_i$ to a parabolic subgroup $P_i\subset Fix(c_i), c_i\in\partial_\infty G/K$, where $G/K$ is the associated symmetric space. (This holds in particular for the inclusion $\rho:\Gamma\rightarrow G$ if 
$int\left(M\right)\cong\Gamma\backslash G/K$ is a ${\bf Q}$-rank $1$ locally symmetric space.) 

Then we have a well-defined simplicial map $B\rho:B\Gamma^{comp}\rightarrow BG^{comp}$ by sending the cone point over $\Gamma_i$ to $c_i$. We will use the shorthand $(B\rho)_*\left[M,\partial M\right]\in H_*^{simp}(BG^{comp})$ for the image of $EM^{-1}\left[M,\partial M\right]\in H_*^{simp}\left(B\Gamma^{comp}\right)$  under $(B\rho)_*$.

\subsection{Proof of Theorem 1}

\begin{pf}
The proof is similar to that of \cite[Theorem 1]{ku2}.
Let $X$ be constructed as in the proof of \hyperref[X]{Lemma \ref*{X}} and let $\partial X=X\cap\left(\partial M\cup\partial M^\tau\right)$. By \hyperref[X]{Lemma \ref*{X}} we have a representation $\rho^X:\pi_1X\rightarrow G$ with $\rho^X\mid_{\pi_1M}=\rho$ and $\rho^X\left(t\right)=A$.
The construction of $X$ implies that
\begin{equation}\label{a}i_{M*}\left[M,\partial M\right]-i_{M^\tau *}\left[M^\tau,\partial M^\tau\right]= i_{T^\tau *}\left[T^\tau,\partial T^\tau\right]\in H_d\left(X,\partial X\right),\end{equation}
where $i_M,i_{M^\tau},i_{T^\tau}$ are the inclusions of $M,M^\tau$ and the mapping torus $T^\tau$ into $X$.

Fix $n\in{\bf N}$ such that $A^n={\bf 1}$ and $\tau^n=id$.
Using the presentation $$\pi_1X=<S,t\mid R, tht^{-1}=\tau_*\left(h\right)\ \forall\ h\in\pi_1\Sigma>$$ we define a surjective homomorphism
$a:\pi_1X\rightarrow {\bf Z}/n{\bf Z}$ by $$a\left(
t\right)=1, a\left(s\right)=0\ \forall\ s\in S.$$ 
Let $\pi:\widehat{X}\rightarrow X$ 
be the $n$-fold cyclic covering with $\Gamma^{\widehat{X}}:=\pi_1\left(\widehat{X},\hat{x}\right)\cong ker\left(a\right)$. 
Consider $\widehat{M}=\pi^{-1}\left(M\right),
\widehat{M}^\tau=\pi^{-1}\left(M^\tau\right)$ and $\Sigma\times{\bf S}^1=\pi^{-1}(T^\tau)$.

The transfer map $tr:H_*\left(X,\partial X\right)
\rightarrow H_*\left(\widehat{X},\partial \widehat{X}\right)$ 
applied to
\hyperref[a]{Equation \ref*{a}} yields
\begin{equation}\label{tr} 
i_{\widehat{M}*}\left[\widehat{M},\partial \widehat{M}\right]-
i_{\widehat{M}^\tau *}\left[\widehat{M}^\tau,\partial \widehat{M}^\tau
\right]= i_{\Sigma\times{\bf S}^1 *}\left[\Sigma\times{\bf S}^1,\partial\Sigma\times{\bf S}^1\right]\in H_d\left(\widehat{X},\partial \widehat{X}\right).\end{equation}

Let $\partial_1M,\ldots,\partial_sM$ be the components of $\partial M$. 
The path $l_i$, which defines the isomorphism of $\pi_1(\partial_iM,x_i)$ with $\Gamma_i\subset\Gamma:=\pi_1(M,x)$ (see Section 3.1) also identifies $\pi_1(\partial_iX)$ to a subgroup $\Gamma_i^X\subset \pi_1(X,x)$ for the corresponding component $\partial_iX$ of $\partial X$.
For each component $\partial_{ik}\widehat{X}$ of $\pi^{-1}(\partial_iX)$ we can use a lift $\hat{l}_{ik}$ of $l_i$ to define an isomorphism from $\pi_1(\partial_{ik}\widehat{X},\hat{x}_{ik})$ to a subgroup
$\Gamma_{ik}^{\widehat{X}}\subset \Gamma^{\widehat{X}}$ 
with $\pi_*\left(\Gamma_{ik}^{\widehat{X}}\right)=\Gamma_i^X$.

The representation $\rho_{\widehat{X}}:=\rho^X\pi_*:\Gamma^{\widehat{X}}\rightarrow
G$ sends all $\Gamma_{ik}^{\widehat{X}}$ to the parabolic group $P_i$. Indeed,  
$\Gamma_i^X=\pi_*\left(\Gamma_{ik}^{\widehat{X}}\right)$ is generated by $t^n$ and elements of $\Gamma_i$, and since 
$\rho_X\left(t^n\right)=A^n={\bf 1}$ we have $\rho_{\widehat{X}}
\left(\Gamma_i^{\widehat{X}}\right)=\rho_X\left(\Gamma_i^X\right)=
\rho_X\left(\Gamma_i\right)=\rho(\Gamma_i)\subset P_i$. Therefore we can
extend $B\rho_{\widehat{X}}$ to 
$$B\rho_{\widehat{X}}:\left(B\Gamma^{\widehat{X}}\right)^{comp}\rightarrow BG^{comp},$$
by mapping the cone point over $B\Gamma_{ik}^{\widehat{X}}$ to $c_i$, the fixed point of the parabolic group $P_i$.

Let $\Psi_{\widehat{X}}:DCone\left(\cup_{i=1}^s\partial_i\widehat{X}\rightarrow \widehat{X}\right)
\rightarrow \mid B\Gamma^{\widehat{X} comp}\mid$ be the extension of the classifying map.
An argument analogous to the closed case (Section 1.4) using $\rho_{\widehat{X}}(\pi_1{\bf S}^1)={\bf 1}$ shows
$$\left(\mid B\rho_{\widehat{X}}\mid\Psi_{\widehat{X}}i_{\Sigma\times{\bf S}^1}
\right)_*=\left(\mid B\rho_{\widehat{X}}\mid\Psi_{\widehat{X}} i_{\Sigma} 
P_1\right)_*:H_d\left(\Sigma\times {\bf S}^1,\partial\Sigma\times{\bf S}^1\right)\rightarrow H_d
\left(\mid BG^{comp}\mid\right).$$
Since the latter homomorphism factors over $H_d\left(\Sigma,\partial\Sigma\right)=0$ we have $$\left(\mid B\rho_{\widehat{X}}\mid
\Psi_{\widehat{X}}i_{\Sigma\times{\bf S}^1}\right)_*\left[\Sigma\times {\bf S}^1,\partial\Sigma\times{\bf S}^1\right]=0$$ 
and
\hyperref[tr]{Equation \ref*{tr}} implies
\begin{equation}\label{d}
\left(\mid B\rho_{\widehat{X}}\mid\Psi_{\widehat{X}}i_{\widehat{M}}\right)_*
\left[\widehat{M},\partial\widehat{M}\right]=\left(
(\mid B\rho_{\widehat{X}}\mid\Psi_{\widehat{X}}
i_{\widehat{M}^\tau}\right)_*\left[\widehat{M}^\tau,\partial \widehat{M}^\tau\right]\in H_d\left(\mid B
G^{comp}\mid\right).\end{equation}

On the other hand we have 
$$(B\rho)_*\left[M,\partial M\right]_{\bf Q}=
\frac{1}{n} (B\rho_X\Psi_Xi_M)_*\pi_*\left[
\widehat{M},\partial\widehat{M}\right]_{\bf Q}=\frac{1}{n}
(B\rho_{\widehat{X}}
\Psi_{\widehat{X}}i_{\widehat{M}})_*
\left[\widehat{M},
\partial \widehat{M}\right]_{\bf Q},
$$
similarly for $M^\tau$, and thus \hyperref[d]{Equation \ref*{d}} implies the claim of Theorem 1.
\end{pf}



\end{document}